\input amstex
\documentstyle{amsppt}
\pagewidth{12.9cm} \pageheight{18.5cm} \topmatter
\title
Trace Formulas on Finite Groups
\endtitle
\author Jae-Hyun Yang
\endauthor
\magnification=\magstep 1 \baselineskip =15mm
\thanks{This work was supported by INHA UNIVERSITY Research
Grant (INHA-22792).
\endgraf
}
\endthanks
\abstract{In this paper, we study the right regular representation
$R_{\Gamma}$ of a finite group $G$ on the vector space consisting
of vector valued functions on $\Gamma\backslash G$ with a subgroup
$\Gamma$ of $G$ and give a trace formula using the work of M.-F.
Vign{\' e}ras. \endgraf Key words and phrases\,:\,representations
of a finite group, trace formula.}
\endabstract
\endtopmatter
\document
\NoBlackBoxes

\def\BC{\Bbb C}

\def\lrt{\longrightarrow}

\def\l{\lambda}

\def\G{\Gamma}

\def\d{\delta}

\def\Gb{\Gamma\backslash G}

\vskip 0.35cm \head {\bf 1.\ Introduction}    
\endhead                
\vskip 0.3cm Following the suggestion of D. Kazhdan, James Arthur
proved the so-called {\it local\ trace\ formula} for a reductive
group $G(F)$ over a non-archimedean local field $F$ investigating
the regular representation of $G(F)\times G(F)$ on the Hilbert
space $L^2(G(F))$\,(cf.\,[1]-[4]\,). Motivated by the work of J.
Arthur on the local trace formula, M.-F. Vign{\' e}ras (cf.\,[10])
gave a trace formula for the regular representation of $G\times G$
in $L^2(G)$ for a finite group $G$. In this paper, motivated by
the above mentioned work of M.-F. Vign{\' e}ras, we study the
trace formula of the right regular representation $R_{\G}$ of a
finite group $G$ on the vector space of consisting of all vector
valued functions on the coset space $\Gb$ for a subgroup $\G.$ We
derive the trace formula for $R_{\G}(f)$ using the result of M.-F.
Vign{\' e}ras (cf.\,[10]). This trace formula simplifies the
proofs of the well known results on a finite group.
\par\medpagebreak
In this paper, we shall study the right regular representation $R$
of $G$ on the vector space $V[\G\backslash G]$ consisting of all
vector valued functions on $\Gb$ with values in $V$ and give a
trace formula for $R_{\G}(f)$ with a function $f$ on $G$. Using
this formula, we derive some well known results.
\par\medpagebreak\noindent
{\bf Notation.} We denote by $\BC$ the complex number field. For a
finite set $A$, we denote by $|A|$ the cardinality of $A$. For a
finite group $G$, we denote by ${\hat G}$ the set of all
equivalence classes of irreducible representations of $G$. For
$\l\in {\hat G}$, we let $d_{\l}$ be the degree of $\l.$

\def\pmn{\par\medpagebreak\noindent}
\par\medpagebreak\noindent

\vskip 0.7cm \head {\bf 2.\ Trace Formula}    
\endhead                 
\vskip 0.3cm Let $\G$ be a subgroup of a finite group $G$. Let $V$
be a finite dimensional complex vector space. We let
$X_{\G}=\G\backslash G$ and denote by $V_{\G}$ the vector space
consisting of all vector valued functions $\varphi:X_{\G}\lrt V$.
We note that $G$ acts on $X_{\G}$ transitively by right
multiplication. We let $R_{\G}$ be the right regular
representation of $G$ on $V_{\G},$ namely,
$$(R_{\G}(g)\varphi)(x)=\varphi(xg),\quad g\in G,\ \varphi\in V_{\G}\ \text{and}\ x\in
X_{\G}.$$ \indent For any $g\in G,$ we set $X_{\G}^g=\{\,x\in
X_{\G}\,|\ xg=x\,\}.$ \pmn {\bf Theorem\ 1.} Let
$G,\,\G,\,V_{\G},\,X_{\G}$ and $X_{\G}^g$ be as above. We let
$\chi_{R_{\G}}$ be the character of the regular representation
$R_{\G}$ of $G$. For each $\l\in {\hat G},$ we let $\chi_{\l}$ be
the character of $\l.$ We assume that $R_{\G}=\sum_{\l\in {\hat
G}}m_{\l}(\G,V)\l$ is the decomposition of $R_{\G}$ into
irreducibles. Here $m_{\l}(\G,V)$ denotes the multiplicity of $\l$
in $R_{\G}$. Then
$$ \chi_{R_{\G}}(g)=\text{dim}_{\BC}\,V\cdot |X_{\G}^g|\quad \text{for\ all}\
g\in G,\tag 1$$
\def\dim{\text{dim}_{\BC}\,V}
\def\ml{m_{\l}(\G,V)}
\def\xg{|X_{\G}^g|}
$$\ml={{\dim}\over {|G|}}\,\sum_{g\in G}\xg\, \chi_{\l}(g^{-1})\quad
\text{for\ each}\ \l\in {\hat G},\tag2$$
\def\G{\Gamma}
\def\lhat{\lambda\in {\hat G}}
$$|G|^2=|\G|\sum_{\lhat}\sum_{g\in G}d_{\l}\,\xg\,
\chi_{\l}(g^{-1}).\tag3$$ \indent For a function $f\in \BC[G],$ we
define the endomorphism $R_{\G}(f)$ of $V_{\G}$ by
$$R_{\G}(f)=\sum_{g\in G}f(g)R_{\G}(g).$$
Then for a function $f\in \BC[G],$
$$\text{tr}\,R_{\G}(f)=\dim\cdot\sum_{g\in G}f(g)\,\xg,\tag4$$
and
for any $f_1,\,f_2\in \BC[G],$
$$\text{tr}\,R_{\{1\}}(f_1*f_2)=\sum_{\lhat}d_{\l}\,\text{tr}\,({\Cal
F}f_1(\l){\Cal F}f_2(\l)).\tag5$$ Here $f_1*f_2$ denotes the
convolution of $f_1$ and $f_2$ defined by
$$(f_1*f_2)(g)=\sum_{h\in G} f_1(h)f_2(h^{-1}g),\qquad g\in G,$$
$d_{\l}$ is the degree of $\l$ and ${\Cal F}f(\l)$ is the Fourier
transform of $f$ defined by
$${\Cal F}f(\l)=\sum_{g\in G}f(g)\,\l(g),\quad \lhat.$$

\pmn {\it Proof.} We let $V[X]$ be the set of all functions
$\phi:X\lrt V$ with values in $V$. We describe a basis for the
vector space $V[X]$ and its dual basis. If $V=\BC$, the vector
space $\BC [X]$ has a basis $\{\,\delta_x\,|\,x\in X\,\},$ where
$$
\delta_x(y):=\cases 1 &\text{if $x=y$}\\
0 &\text{otherwise}.\endcases
$$
For $x\in X$ and $v\in V$, we define the function $\delta_x\otimes
v:X\lrt V$ by
$$
(\d_x\otimes v)(y):=\cases v &\text{if $x=y$}\\
0 &\text{otherwise}.\endcases
$$
Let $\{\,v_1,\cdots,v_n\,\}$ be a basis for $V$ with
$\text{dim}_{\BC}V=n.$ Then it is easy to see that the set
$\{\,\d_x\otimes v_k\,|\,x\in X,\ 1\leq k\leq n\,\}$ forms a basis
for $V[X]$. Let $V^*$ be the dual space of $V$. For $x\in X$ and
$v^*\in V^*,$ we define the linear functional $\d_x^*\otimes
v^*:V[X]\lrt \BC$
$$\left( \d_x^*\otimes v^*\right)(\phi):=<\phi(x),\,v^*>,\quad \phi\in
V[X].$$ Suppose $\{\,v_1^*,\cdots,v_n^*\,\}$ is the dual basis of
a basis $\{\,v_1,\cdots,v_n\,\}$. Then we see easily that the set
$\{\,\d_x^*\otimes v_k^*\,\vert\ x\in X,\ 1\leq k\leq n\,\}$ forms
a basis for the dual space $V[X]^*$. We also see that for each
$g\in G,$
$$<R_{\G}(g)(\d_x\otimes v_k),\,(\d_x^*\otimes v_k^*)>=\cases 1 &
\text{if $xg=x$}\\
0 &\text{otherwise}.\endcases$$ Therefore
$\chi_{R_{\G}}(g)=\text{tr}\,R_{\G}(g)=n\cdot \xg$ for each $g\in
G$. This proves Formula (1).
\par\smallpagebreak
We define the hermitian inner product $<\ ,\ >$ on the group
algebra $\BC[G]$ by
$$<f_1,\,f_2>={1\over {|G|}}\,\sum_{g\in G}f_1(g)\,{\overline
{f_2(g)}},\quad f_1,\,f_2\in \BC[G].$$ Since
$\chi_{R_{\G}}=\sum_{\lhat}\ml\,\chi_{\l},$ we have
$$\align
\ml&=<\chi_{R_{\G}},\chi_{\l}> \quad (\,\text{by\ Schur orthogonality relation\,(cf.\,[6],\,p.\,148)})\\
&={1\over {|G|}}\sum_{g\in
G}\chi_{R_{\G}}(g)\,\chi_{\l}(g^{-1})\\
&={n\over {|G|}}\sum_{g\in G}\xg\,\chi_{\l}(g^{-1}).\qquad\quad
(\,\text{by (1)}\,)
\endalign$$
This proves Formula (2).
\par\smallpagebreak
We observe that
$$\text{dim}_{\BC}\,V_{\G}={{|G|}\over {|\G|}}\cdot \dim
\tag6$$ Since $R_{\G}=\sum_{\l\in {\hat G}}m_{\l}(\G,V)\l$, we
see that
$$\text{dim}_{\BC}\,V_{\G}=\sum_{\lhat}\ml\cdot d_{\l},\tag7$$
where $d_{\l}$ denotes the degree of $\lhat$. By substituting (2)
into (7), we get
$$\text{dim}_{\BC}\,V_{\G}={{\dim}\over {|G|}}\cdot \sum_{\lhat}\sum_{g\in G}
d_{\l}\,\xg\,\chi_{\l}(g^{-1}).\tag8$$ Therefore according to (6)
and (8), we obtain Formula (3).
\par\smallpagebreak
Let $f\in \BC[G].$ Then we obtain
$$\align \text{tr}\,R_{\G}(f)&=\text{tr}\,\left(\sum_{g\in
G}f(g)R_{\G}(g)\right)\\
&=\sum_{g\in G}f(g)\,\text{tr}\,(R_{\G}(g))\\
&=\sum_{g\in G}f(g)\chi_{R_{\G}}(g)\\
&=\dim\cdot\sum_{g\in G}f(g)\,\xg.\qquad (\,\text{by}\ (1)\,)
\endalign$$
This proves Formula (4).
\par\medpagebreak
Finally we shall prove Formula (5). If we take $\G=\{1\},$ then
$X_{\G}^1=G$ and $X_{\G}^g=\emptyset$ if $g\neq 1.$ Thus by
Formula (4), we get
$$\text{tr}\,R_{\{1\}}(f)=\dim\cdot |G|\,f(1).$$
Therefore
$$f(1)={{\text{tr}\,R_{\{1\}}(f)}\over {|G|\,\dim}}.\tag9$$
We recall the fact\,(\,see [6], Corollary 3.4.5\,) that for any
$f_1,\,f_2\in \BC[G]$, the following Plancherel formula holds\,:
$$(f_1*f_2)(1)=
{{\text{tr}\,R(f_1*f_2)}\over {|G|}}={1\over {|G|}}
\sum_{\lhat}d_{\l}\,\text{tr}_{V_{\l}}\,({\Cal F}f_1(\l){\Cal
F}f_2(\l)),\tag10$$ where $\text{tr}_{V_{\l}}(A)$ denotes the
trace of an endomorphism $A:V_{\l}\lrt V_{\l}$ with the
representation space of $\l$.
\par\smallpagebreak
On the other hand, for for any $f_1,\,f_2\in \BC[G]$, we get
$$(f_1*f_2)(1)={{\text{tr}\,R_{\{1\}}(f_1*f_2)}\over {|G|}}\qquad
(\,\text{by}\ (9)\,).\tag11$$ Hence according to (10) and (11), we
obtain Formula (5). \hfill $\square$
\par\bigpagebreak\noindent
{\bf Corollary\ 2.} (a)  $|G|=\sum_{\lhat}d_{\l}^2.$ \pmn (b)
$|G|=\sum_{g\in G}\sum_{\lhat}\,d_{\l}\,\chi_{\l}(g^{-1}).$  \pmn
(c) Let $\l_0$ be the trivial representation of $G$. Then
$$m_{\l_0}(\G,V)={{\dim}\over {|G|}}\cdot \sum_{g\in G}\xg.$$
In particular, $m_{\l_0}(\{1\},\BC)=1.$ \pmn (d) For any $f\in
\BC[G]$,
$$f(1)={{\text{tr}\,R_{\{1\}}(f)}\over {|G|\cdot \dim}}.$$
(e) For a subgroup $\G$ of $G$,
$$\sum_{\lhat} \ml^2={{({\dim})^2}\over {|G|}}\cdot \sum_{g\in G}
\xg\,|X_{\G}^{g^{-1}}|.$$  (f) For each $\lhat,\
m_{\l}(\{1\},V)=d_{\l}\, \dim\neq 0.$ That is, each $\lhat$ occurs
in the regular representation $(R_{\{1\}},\,V[G])$ of $G$ with
multiplicity $d_{\l}\cdot\dim$.  \pmn (g) For each $\lhat$,
$$\sum_{g\in G} \chi_{\l}(g)={{|G|}\over {\dim}}\cdot
m_{\l}(G,V).$$ {\it Proof.} (a) If we take $\G=\{1\},$ we see
easily that $X_{\{1\}}^1=G$ and $X_{\{1\}}^g=\emptyset$ if $g\neq
1.$ Then we get
$$\align
|G|^2&=1\cdot \sum_{\l\in {\hat G}}d_{\l}\cdot |X_{\{1\}}^1|\cdot
\chi_{\l}(1)\qquad (\,\text{by}\ (3)\,)\\
&=|G|\,\sum_{\lhat}d_{\l}^2. \qquad (\,\text{because}\
\chi_{\l}(1)=d_{\l}\,)
\endalign$$
This proves Formula (a). We recall that another proof of (a)
follows from the fact that the group algebra $\BC[G]$ is
isomorphic to $\sum_{\lhat}\text{End}(V_{\l})$ as algebras, where
$V_{\l}$ is the representation space of $\lhat$ (cf. [9]). \pmn
(b) We take $\G=G.$ It is easy to see that $X_G=\{{\bar 1}\}$ is a
point and $X_G^g=X_G$ for all $g\in G.$ According to Formula (3),
we obtain
$$|G|^2=|G|\,\sum_{\lhat}\sum_{g\in G}d_{\l}\chi_{\l}(g^{-1}).$$
This proves the statement (b). \pmn (c) It follows from Formula
(2). \pmn (d) We take $\G=\{1\}.$ Then $X_{\G}^1=G$ and
$X_{\G}^g=\emptyset$ for $g\neq 1.$ From Formula (4), we obtain
$$\text{tr}\,R_{\{1\}}(f)=\dim\cdot |G|\,f(1).$$
\pmn (e) By Schur orthogonality relation,
$<\chi_{R_{\G}},\,\chi_{R_{\G}}>=\sum_{\lhat}\ml^2.$ On the other
hand, according to the formula (1),
$$<\chi_{R_{\G}},\,\chi_{R_{\G}}>={1\over {|G|}}\,\sum_{g\in
G}\chi_{R_{\G}}(g)\chi_{R_{\G}}(g^{-1})={{({\dim})^2}\over
{|G|}}\cdot \sum_{g\in G} \xg\,|X_{\G}^{g^{-1}}|.$$ \pmn (f) We
take $\G=\{1\}.$ Then $X_{\{1\}}^1=G,\ X_{\{1\}}^g=\emptyset$ for
$g\neq 1$, and $V_{\{1\}}=V[G]$. Therefore we obtain the desired
result from Formula (2).

\pmn (g) We take $\G=G.$ We see easily that $|X_G^g|=1$ for all
$g\in G.$ According to Formula (2), we get
$$m_{\l}(G,V)={{\dim}\over {|G|}}\sum_{g\in G}\chi_{\l}(g^{-1})=
{{\dim}\over {|G|}} \sum_{g\in G}\chi_{\l}(g).$$ Hence this proves
the statement (g).
 \hfill $\square$
\par\bigpagebreak\noindent
{\bf Theorem\ 3.} For each $f\in \BC[G],$ we have the following
trace formula
$$\text{tr}\,R_{\G}(f)={{\dim}\over {|G|}}\cdot \sum_{g\in
G}\xg\,|Z_g|\,f(C_g),\tag12$$ where $Z_g$ is the centralizer of
$g$ in $G,\ C_g$ is the conjugacy class of $g$ and
$f(C_g)=\sum_{h\in C_g}f(h).$

\pmn {\it Proof.} For $f\in \BC[G]$ and $\lhat$, we define
$$\l (f):=\sum_{g\in G} f(g)\l(g).$$
\indent Investigating the spectral decomposition of the regular
representation $R$ of $G\times G$ on $\BC [G]$ defined by
$$\left( R(g_1,g_2)F\right)(g)=F(g_1^{-1}gg_2),\quad g,g_1,g_2\in
G,\ F\in \BC[G].$$ M.-F. Vign{\'e}ras (cf.\,[10], p.284) obtained
the following trace formula
$$
|Z_g|F(C_g)=\sum_{\pi\in
\hat{G}}\chi_{\pi}(g^{-1})\,\text{tr}\,\pi(F)\tag 13$$ for any
$g\in G$ and $F\in \BC[G].$

Let $R_{\G}=\sum_{\l\in {\hat G}}m_{\l}(\G,V)\l$ be the
decomposition of $R_{\G}$ into irreducibles. If $f\in \BC[G],$
$$\align
\text{tr}\,R_{\G}(f)&=\sum_{\lhat}\ml\,\text{tr}\,\l(f)\\
&={{\dim}\over {|G|}}\sum_{\lhat}\sum_{g\in
G}\xg\,\chi_{\l}(g^{-1})\,\text{tr}\,\l(f) \qquad (\,\text{by (2)}\,) \\
&={{\dim}\over {|G|}}\,\sum_{g\in G}\xg\,\left(\sum_{\lhat}
\chi_{\l}(g^{-1})\,\text{tr}\,\l(f)\right)\\
&={{\dim}\over {|G|}}\cdot \sum_{g\in
G}\xg\,|Z_g|\,f(C_g).
\endalign$$
The last equality follows from Formula (13). \hfill $\square$
\par\bigpagebreak\noindent
{\bf Corollary\ 4.} Let $\G$ be a subgroup of $G$. Then for any
$f\in \BC[G]$, we have the following identity
$$ |G|\sum_{g\in G}f(g)|X_{\G}^g|=\sum_{g\in G}|X_{\G}^g||Z_g|
f(C_g),\tag13$$ where $Z_g$ is the centralizer of $g$ in $G,\ C_g$
is the conjugacy class of $g$ and $f(C_g)=\sum_{h\in C_g}f(h).$
\pmn {\it Proof.} The proof follows immediately from Formula (4)
and the trace formula (12). \hfill $\square$

\pmn {\bf Remark\ 5.} The trace formula (12) is similar to the
trace formula on the adele group. For the trace formula on the
adele group, we refer to [1]-[5],\,[7] and [8].

\par\medpagebreak\noindent
{\bf Remark\ 6.} If $\G\neq \{1\},$ the multiplicity $\ml$ of some
$\lhat$ may be zero. It is natural to ask when $\ml$ is not zero.
Namely, which $\lhat$ does occur in the regular representation
$(R_{\G},\,V_{\G})$ of $G$\,?

\vskip 0.5cm

\par\bigpagebreak\noindent
Department of Mathematics\par\noindent Inha
University\par\noindent Incheon 402-751\par\noindent Republic of
Korea \vskip 0.25cm \noindent {\tt email\
address\,:\,jhyang\@inha.ac.kr}


\vskip 1cm \Refs \widestnumber\key{\bf C-R}

\ref\key{\bf 1} \by J. Arthur \paper The trace formula and Hecke
operators: Number theory, trace formulas and discrete groups
(Oslo, 1987)\publ Academic Press
 \yr 1989\pages 11-27
\endref

\ref\key{\bf 2} \bysame \paper Towards a local trace formula:
algebraic analysis, geometry and number theory (Baltimore, MD,
1988) \publ Johns Hopkins Univ. Press, Baltimore, MD \yr 1989
\pages 1-23
\endref

\ref\key{\bf 3} \bysame \paper Some problems in local harmonic
analysis : harmonic Analysis on reductive groups (Brunswick, ME,
1989) \publ Prog. Math. Birkh{\"a}user, Boston vol 101 \yr 1991
\pages 57-78
\endref


\ref\key{\bf 4} \bysame \paper A local trace formula \publ IHES
Publication Math.\vol 73 \yr 1991 \pages 5-96
\endref






\ref\key{\bf 5} \by S. Gelbart \book Lectures on Arthur-Selberg
Trace Formula \publ American Math. Soc., Providence \yr 1996
\endref



\ref\key{\bf 6} \by R. Goodman and N. Wallach \book
Representations and Invariants of the Classical Groups\publ
Cambridge University Press  \yr 1998
\endref


\ref\key{\bf 7} \by A. W. Knapp \paper Theoretical Aspect of the
Trace Formula for $GL(2)$, Proc. Symosia in Pure Math. \publ
American Math. Soc. {\bf 61} \yr 1997\pages 355-405
\endref

\ref\key{\bf 8} \by A. W. Knapp and J. D. Rogawski\paper
Applications of the Trace Formula, Proc. Symosia in Pure Math.
\publ American Math. Soc. {\bf 61} \yr 1997\pages 413-431
\endref




\ref\key{\bf 9} \by J.-P. Serre \book Linear representations of
finite groups \publ Springer-Verlag \yr 1977
\endref
\ref\key{\bf 10} \by M.-F. Vign{\'e}ras \paper An elementary
introduction to the local trace formula of J. Arthur. The case of
finite groups  \jour Jber. d. Dt. Math.-Verein.,
Jubil{\"a}umstagung 1990, B.G. Teubner Stuttgart  \pages 281-296
\yr 1992
\endref
\endRefs
\end{document}